\documentclass{article}
\usepackage{amsfonts}
\usepackage{amssymb}

\newtheorem{theorem}{Theorem}

\newtheorem{definition}[theorem]{Definition}
\newtheorem{example}[theorem]{Example}

\newtheorem{proposition}[theorem]{Proposition}
\newtheorem{remark}[theorem]{Remark}

\input{tcilatex}
\begin{document}

\title{\textbf{A new approach for Imprecise Probabilities}}
\author{Marcello Basili \\
DEPS University of Siena\\
marcello,basili@unisi.it \and Luca Pratelli \\
Naval Academy Leghorn\\
luca\_pratelli@marina.difesa.it}
\maketitle

\begin{abstract}
This paper introduces a novel concept of interval probability measures that enables the representation of imprecise probabilities, or uncertainty, in a natural and coherent manner. Within an algebra of sets, we introduce a notion of weak complementation denoted as $\psi$. The interval probability measure of an event $H$ is defined with respect to the set of indecisive eventualities $(\psi(H))^c$, which is included in the standard complement $H^c$.

We characterize a broad class of interval probability measures and define their properties. Additionally, we establish an updating rule with respect to $H$, incorporating concepts of statistical independence and dependence. The interval distribution of a random variable is formulated, and a corresponding definition of stochastic dominance between two random variables is introduced. As a byproduct, a formal solution to the century-old Keynes-Ramsey controversy is presented.
\end{abstract}

\bigskip\section{Introduzione}

A highly sophisticated literature is dedicated to addressing the representation of uncertainty or ambiguity as quite different from risk represented in the standard probability theory. Uncertain representation encompasses various approaches such as belief functions and plausibility functions (Smets 1981, Dempster 1967, Shafer 1976), upper and lower probabilities (Boole 1958,[1854]), Kyburg 1961, Fagin et al. 1995), and sets of probability measures (Good 1962, Gilboa and Schmeidler 1989).

In the context of additive probability measures, the Boolean Algebra provides the mathematical framework
where the properties of Ramsey's probability concept are formally defined.
The properties of additivity and unique complementation, inherent to Boolean
Algebra, give rise to the Law of Excluded Middle (Ramsey 1926, 65).

In 1936, Stone introduced a structural analysis of Boolean algebras as an
organic system, an alternative to studies of symbolic logic and \lq\lq their
preoccupation with algorithms, identities, and equations, or with the
logical interrelations of the formal properties of the various Boolean
operations''(Stone 1938, 807). In his Theorem of Representation, Stone
demonstrated that every abstract Boolean algebra or complemented
distributive and bounded lattice corresponds isomorphically to some algebra
of sets with the usual set-theoretic operations.

During the same year, amidst the Einstein-Podolsky-Rosen and Bohr's
controversy regarding the quantum mechanical description of physical
reality, Birk\-hoff and von Neumann published "The Logic of Quantum
Mechanics." They formulated a probability theory for quantum mechanics on
non-distributive and orthocomplemented lattices, where the property of
additivity is not satisfied.

Koopman emancipated the notion of probability from epistemology and
simultaneously provided an axiomatic foundation that addresses the
limitation found in "the fact that all the axiomatic treatments of intuitive
probability current in the literature take as their starting point a number
(usually between 0 and 1) corresponding to the degree of rational belief or
credibility of the event in question" (Koopman 1940, 269). In this
perspective, Koopman defined numerical probability and its properties on a
Boolean algebra.

Apart from isomorphisms, these theories, which incorporate viewpoints from Ramsey, De Finetti, Savage and others, share a fundamental starting point focused on
the consideration of complementation - the usual
set-theoretic complementation - whithin an algebra of sets utilizing the
operations of union and intersection.

This paper introduces weak complementation, a distinct notion of set- theoretic complementation, in the algebra of sets (or events) of a generic non-empty set. Under this definition, the algebra of sets transforms into a weakly complemented lattice that deviates from De Morgan's rules and lacks (unique) complementation. Nevertheless, this framework facilitates the conceptualization of uncertainty associated with an event in a natural and coherent manner. 'Imprecise belief' is thereby represented by the measure of the event characterizing the uncertainty, providing the opportunity to define general classes of interval probability measures.

By interpreting eventualities as the causal factors leading to the occurrence of a given random phenomenon, we introduce the set of indecisive eventualities (uncertain opportunities) associated with an event $H$. This set constitutes the most comprehensive collection of factors, not necessarily void, that could contribute to the realization of the random phenomenon when considering the causes determining $H$. Such an indecisive set provides an intuitive definition of uncertainty in a weakly complemented distributive algebra ${\mathcal A}$ of events.

The concept of weak complementation 
also allows the transformation of ${\mathcal A}$ into a complemented, non-distributive lattice by modifying only the partial ordering relation $\subseteq$ and the standard set union $\cup$. In this context, the empty set $\emptyset$ does not yet serve as the minimum according to the modified partial order relation. The uncertainty of an event $H$ can be represented by eventualities of $H\cup^\prime\emptyset$ that do not belong to $H$.

Unlike other approaches, we establish uncertainty in both a distributive algebra that is only weakly complemented and, equivalently, in a non-distributive complemented algebra. Importantly, these two approaches lead to the same notion of interval probability measure for $H$ since it is defined in relation to its set of indecisive eventualities. As this indecisive set expands, the width of the interval probability measure also increases. The interval probability of $H$ (or imprecise probability of $H$) according to $P$ represents the confidence interval for the occurrence of a random phenomenon, considering that the causes determining $H$ and its set of indecisive eventualities are measured by $P$.

In this perspective, two more general definitions of interval (probability)
measures are provided, considering a random variable that generalizes the
set of uncertain opportunities and a capacity. These new concepts can be
readily applied to real-world problems involving imprecise probabilities,
such as the IPCC's Reports.

Furthermore, we characterize specific classes of interval probability
measures and, in the context of multiple interval probability measures (collections), we
establish an interval probability measure for the collection, defined
alongside the measure of coherence. This introduces an intuitive method for
aggregating and evaluating uncertain probability judgments.

 Additionally, we introduce fundamental concepts—conditional interval probability measures and stochastic independence of an event $A$ conditioned on $H$. In contrast to methods that rely on assumptions such as rectangularity, menu dependence, change of subjective perception, and others that are reducible to a Boolean algebra (as in Dempster-Shaffer's rule), our novel updating rule for determining the probability of $A$ conditioned on $H$ is presented without any algebraic or probabilistic assumptions. It is represented by the interval with extremes derived through the standard conditional probability of $A$ and $A \cup \tilde A$ concerning $H \cup \tilde H$, where $\tilde A$ and $\tilde H$ denote the events of indecisive eventualities for $A$ and $H$, respectively. The
notion of independence is thereby established (e.g. Kolmogorov), and it is not
necessarily symmetric; in other words, $A$ can be independent of $H$, but $H$
is not necessarily independent of $A$.

\smallskip As a consequence of the concept of the interval distribution of a
random variable, we introduce a novel notion of stochastic dominance between
two random variables based on interval probability measures. 

Finally, the representation of uncertainty by sets and the new interval
probability measure are used to provide a formal resolution to the
historical controversy between Keynes and Ramsey.

The organization of the paper unfolds as follows: Section 2 establishes
crucial aspects in the class of algebraic structures and defines the concept
of weak relative complementation. Interval probability measures are
characterized and discussed in Section 3 while Section 4 and Sections 5  introduce conditional interval probability
measures and the concept of conditional independence. In Section 6, the interval distribution of a random variable is set up and a new notion of stochastic dominance between random variables is introduced. Theoretical properties
of these \textsl{imprecise probability measures}, along with illustrative
examples, are explored in each respective section. Section 7 formulates a possible resolution of the Keynes-Ramsey's century-hold controversy, and Section 8 presents some conclusions. Any results are proven in the
Appendix.

\medskip

\section{Weak and ${\mathcal Z}$-complementation}

Let $\Omega$ be a generic non-empty set, and consider an algebra $\mathcal{A}
$, sometimes a $\sigma$-algebra, of subsets of $\Omega$, denoted as events. Equipped with the relation of inclusion (order relation) and the standard operations of intersection and union, $\mathcal{A}$ constitutes a bounded distributive lattice. For any pair of events $H$ and $K$, let $H^c$ denote the usual set-theoretic complement and $K\setminus H = K \cap H^c$.

\begin{definition}
A map $\psi:{\mathcal A}\mapsto{\mathcal A}$ is said to be a weak complementation if the following conditions are satisfied:

\smallskip $(a)$ \ $\psi(H)\subseteq H^c$ for any event $H\in{\mathcal A}$,

\smallskip $(b)$  For any pair of events $H,K\in{\mathcal A}$, with $H\subseteq K$, it holds 
$$
K_\psi\subseteq H_\psi
$$where $K_\psi=K^c\setminus\psi(K)$ and $H_\psi=H^c\setminus\psi(H)$.

\medskip\noindent Moreover, $\psi$ is a regular weak complementation if $\psi$ also satisfies the condition 

\medskip $(c)$ $\psi\big(\psi(H)\big)\subseteq H,$ for any event $H$. 
\end{definition}

\medskip\noindent The event $H_{\psi} = H^c \setminus \psi(H)$ is referred to as the \textit{set of uncertain opportunities (or indecisive eventualities) of $H$ relative to weak complementation $\psi$}. For each event $H$, there exists a partition of $\Omega$ given by 
\begin{equation}
\Omega = H \cup \psi(H) \cup H_{\psi}, \end{equation}
where $H_{\psi}$ is not necessarily the empty set. If eventualities are interpreted as the causes leading to the occurrence of
a given random phenomenon, then the set $H \cup H_{\psi}$ is the
largest event that can realize the random phenomenon when the causes
determine $H$.

\smallskip In a classical setting, $\psi(H)$ coincides with $H^c$, and the set of uncertain opportunities is always the empty set. However, in a regular weak complementation setting, if $H \subseteq K$, it implies $H\cup\psi(H)\subseteq K\cup\psi(K)$; nevertheless, it does not guarantee $\psi(K) \subseteq \psi(H)$, and De Morgan's rules do not necessarily hold. 

\smallskip
Now, a significant example of regular weak complementation is proposed. To this aim, let's consider a finite family $\mathcal{Z} = \{Z_0, \ldots, Z_M\}$ of non-empty events that are pairwise disjoint, with $\bigcup_{n=0}^M Z_n = \Omega$, where $M$ is a natural number. More precisely,

\begin{definition} The map $\psi_{\mathcal Z}$ defined by 
\begin{equation}
\psi_{\mathcal Z}(H) = \bigcup_{n : H \cap Z_n \neq \emptyset} H^c \cap Z_n 
\end{equation}
is said to be the  $\mathcal{Z}$-complementation. Moreover, $\psi_{\mathcal Z}(H)$ is the complementary event of $H$ with respect to $\mathcal{Z}$ and the $\mathcal Z$-uncertainty of $H$ is determined by 
\begin{equation}
H_{\psi_{\mathcal{Z}}}=\bigcup_{n : H \cap Z_n=\emptyset}\,  Z_n.
\end{equation}
\end{definition}

\smallskip\noindent It is easy to observe the set of uncertain opportunities $H_{\psi_{\mathcal{Z}}}$ coincides with the event $\cup_{n:Z_n\subseteq H^c}\, Z_n$ and, for every $n=0,\ldots,M$ it holds
\begin{equation}
\psi_{\mathcal Z}(Z_n)=\psi_{\mathcal Z}(\emptyset) =\psi_{\mathcal{Z}}(\Omega)=\emptyset.
\end{equation}  

\medskip\begin{remark}
The implication $H\subset K \Longrightarrow H \cup H_{\psi_{\mathcal{Z}}} \subseteq K \cup K_{\psi_{\mathcal{Z}}}$ is not true if $M \geq 1$. To obtain a counterexample, take a
non-empty event $B$ strictly contained in $Z_1$ and consider $H = Z_0$ and $K = Z_0\cup B$. In this case, $\psi_{\mathcal Z}(H) = \emptyset$ and $\psi_{\mathcal Z}(K) = Z_1 \setminus B$, violating not only the relation $H \cup H_{\psi_{\mathcal{Z}}} \subseteq K \cup K_{\psi_{\mathcal{Z}}}$ but also $\psi_{\mathcal Z}(K) \subseteq \psi_{\mathcal Z}(H)$. However, indecisiveness narrows down, i.e.,
$K_{\psi_{\mathcal{Z}}}\subseteq H_{\psi_{\mathcal{Z}}}$ and it holds $H\cup \psi_{\mathcal{Z}}(H)\subseteq K\cup \psi_{\mathcal{Z}}(K)$
\end{remark}

\smallskip Despite these
considerations, De Morgan's rules are not satisfied, as evidenced by the following result:

\begin{proposition}
$\psi_{\mathcal Z}$ is a regular weak complementation and, for any pair $H,K$ of events, the following identities are true:

\smallskip\noindent

\begin{equation}
\psi_{\mathcal Z}(H\cup K)=H^c\cap K^c\cap \bigcup_{n: H\cap Z_n \neq
\emptyset\, \vee\, K\cap Z_n \neq \emptyset}\, Z_n
\end{equation}

\begin{equation} \psi_{\mathcal Z}(H)\cap\psi_{\mathcal Z}(K)=H^c\cap K^c\cap \bigcup_{n: H\cap
Z_n \neq \emptyset\, \wedge\, K\cap Z_n \neq \emptyset}\, Z_n
\end{equation} 
Moreover, 
\begin{equation}
\psi_{\mathcal Z}(H\cap K)= \bigcup_{n: H\cap K\cap Z_n \neq \emptyset}\,
(H^c\cup K^c)\cap Z_n 
\end{equation}
and 
\begin{equation}\psi_{\mathcal Z}(H)\cup \psi_{\mathcal Z}(K) =J\cup \bigcup_{n: H\cap Z_n \neq
\emptyset\, \wedge\, K\cap Z_n \neq \emptyset}\, (H^c\cup K^c)\cap Z_n
\end{equation}
where $
J=\cup_{n: H\cap Z_n \neq \emptyset=K\cap Z_n}\, H^c\cap Z_n \, \cup_{n:
K\cap Z_n \neq \emptyset=H\cap Z_n}\, K^c\cap Z_n.$

\medskip\noindent In particular
\begin{equation}
\psi_{\mathcal Z}(H)\cap \psi_{\mathcal Z}(K)\subseteq \psi_{\mathcal Z}(H\cup K),\qquad \psi_{\mathcal Z}(H\cap K)\subseteq \psi_{\mathcal Z}(H)\cup \psi_{\mathcal Z}(K).
\end{equation}
\end{proposition}

\medskip\noindent It is important to note that $H_{\psi_{\mathcal{Z}}}$ is the
empty set if and only if the cardinality of $\{n: H \cap Z_n \neq
\emptyset\} $ is $M+1$. Indeed, if (and only if) the cardinality is $M+1$,
it holds $\cup_{n: H \cap Z_n \neq \emptyset} Z_n = \Omega$. Roughly speaking, the set ${\psi_{\mathcal Z}}(H)$ becomes closer to the set-theoretic complement $H^c$ as the cardinality of $\{n: H \cap Z_n \neq \emptyset\}$ increases. 
An interesting context, both from a theoretical standpoint and for applications, is already established for $M=1$. Let's explore how the concept of weak complementation operates in this scenario.

\medskip\begin{example}
Explore the sample space $\Omega =\{ {\omega_{01}, \omega_{10}, \omega_{00}, \omega_{11}}\}$, where each $\omega_{ij}$, with $i, j = 0,1$, represents an eventuality interpreted as either the causal factors $i$ and $j$ or the outcomes $i$ and $j$ obtained in a non-exact manner. These events contribute to the determination of a given random phenomenon.

However, due to the imprecision, indecisive eventualities arise when considering an event $H$. For instance, if $H=\{\omega_{ij}\}$, there exists a unique complementary eventuality $\omega_{i^\prime j^\prime}$, where $|i-i^\prime|+|j-j^\prime|=2$. Meanwhile, eventualities $\omega_{i^{\prime\prime}j^{\prime\prime}}$ with $|i-i^{\prime\prime}|+|j-j^{\prime\prime}|=1$ are deemed indecisive (or uncertain opportunities) concerning $\omega_{ij}$. In simpler terms, to negate $\omega_{ij}$, one needs to modify both indices $i$ and $j$. 
Within this framework, a partition ${\mathcal Z}=\{Z_0,Z_1\}$ is established, where
$$Z_0=\{\omega_{01},\omega_{10}\},\qquad Z_1=\{\omega_{11},\omega_{00}\}$$ Consequently, a weak complementation is defined by $\phi_{\mathcal{Z}}$, precisely as in $(2)$. Specifically, the set of indecisive eventualities within $\{\omega_{01}\}$ corresponds to the event $Z_1$ while $\{\omega_{01},\omega_{11}\}_{\psi_{\mathcal Z}}=\{\omega_{01},\omega_{00}\}_{\psi_{\mathcal Z}}=\emptyset$.  These conclusions will be significant in Section 8 when we delve into  the \lq barometer-dark clouds setting\rq. To obtain the conventional complementation and the classic interpretation, satisfying the Law of Excluded Middle, one should clearly consider the degenerate partition  ${\mathcal Z}^\prime=\{\Omega\}$.

\end{example}

\smallskip  It is worth noting that if ${\mathcal{Z}}^\prime$ is a finer partition of $\Omega$ than ${\mathcal{Z}}$, the complementary event $\psi_{{\mathcal{Z}}^\prime}(H)$ is contained in $\psi_{{\mathcal{Z}}}(H)$ for any event $H$. Consequently, $H_{\psi_{\mathcal{Z}}} \subseteq H_{\psi_{\mathcal{Z}}^\prime}$, leading to an increase in indecisiveness.

\smallskip Moreover, if ${\mathcal{A}}$ is a $\sigma$-algebra of subsets of $\Omega$, ${\mathcal{Z}}$ can be a partition, not necessarily countable, of events in $\Omega$, as is the case when $\mathcal{A}$ is the class of all subsets of $\Omega$. Furthermore, the definition could be
extended to a class ${\mathcal{Z}}$ of events that are not necessarily
pairwise disjoint. In this case, $(2)$ remains unchanged, and  $\psi_{\mathcal Z}$ is a weak but not regular complementation.

\medskip Finally, it is crucial to observe how weak regular complementation allows the transformation of ${\mathcal A}$ into a complemented non-distributive lattice by modifying only the partial ordering relation $\subseteq$ and the usual set union $\cup$. More precisely, given $\psi_{\mathcal Z}$, consider the relation $\subseteq_{\psi_{\mathcal Z}}$ and the two set operations $\cup_{\psi_{\mathcal Z}},\cap_{\psi_{\mathcal Z}}$ defined by
\begin{equation}
H\subseteq_{\psi_{\mathcal Z}}K\Longleftrightarrow      {\psi_{\mathcal Z}}(K)\subseteq {\psi_{\mathcal Z}}(H), \ \ H\cup_{\psi_{\mathcal Z}} K=[\psi_{\mathcal Z}(H\cup K)]^c,\ \ \cap_{\psi_{\mathcal Z}} =\cap.
\end{equation}
As noted in Remark 4,  the partial ordering relations $\subseteq_{\psi_{\mathcal Z}}$ and $\subseteq$ are two distinct relations that do not imply each other. The set ${\mathcal A}$, equipped with $\subseteq_{\psi_{\mathcal Z}}$, $\cup_{\psi_{\mathcal Z}}$, and $\cap_{\psi_{\mathcal Z}}$, forms a complemented lattice with complementation provided by $\psi_{\mathcal Z}$ which satisfies
\begin{equation}
H\cup_{\psi_{\mathcal Z}}\psi_{\mathcal Z}(H)=\Omega, \quad H\cap_{\psi_{\mathcal Z}}\psi_{\mathcal Z}(H)=\emptyset, \quad \psi_{\mathcal Z}(\psi_{\mathcal Z}(H))\cap_{\psi_{\mathcal Z}}H^c=\emptyset
\end{equation}
Moreover, the empty set is not the minimum of ${\mathcal A}$ concerning $\subseteq_{\psi_{\mathcal Z}}$, while $\Omega$ serves as the maximum.  Despite the operations $\cup_{\psi_{\mathcal Z}},\cap_{\psi_{\mathcal Z}}$ satisfying the properties of associativity, and for any event $H$, it follows that 
\begin{equation}
H\cup_{\psi_{\mathcal Z}} H=H\cup_{\psi_{\mathcal Z}} \emptyset=
H\cup H_{\psi_{\mathcal Z}},
\end{equation}  the complemented lattice $({\mathcal A},\cup_{\psi_{\mathcal Z}},\cap_{\psi_{\mathcal Z}},\subseteq_{\psi_{\mathcal Z}},{\psi_{\mathcal Z}})$ is a non-distributive  algebra of subsets since $H\cap_{\psi_{\mathcal Z}}(H\cup_{\psi_{\mathcal Z}}\psi_{\mathcal Z}(H))\neq (H\cap_{\psi_{\mathcal Z}}H)\cup_{\psi_{\mathcal Z}}(H\cap_{\psi_{\mathcal Z}}\psi_{\mathcal Z}(H)).$

\medskip Furthermore, it shares many similarities with, but involves fewer axioms (such as irreducibility) compared to the non-distributive and orthocomplemented lattice introduced in 'The Logic of Quantum Mechanics' by Birkhoff and Von Neumann. Thus, the definition of interval probability measure (introduced in the next Section) remains consistent whether considering $({\mathcal A},\subset,\cup,\cap)$ with weak complementation $\psi$ or the non-distributive algebra $({\mathcal A},\cup_{\psi_{\mathcal Z}},\cap_{\psi_{\mathcal Z}},\subseteq_{\psi_{\mathcal Z}},{\psi_{\mathcal Z}})$.

\section{Interval Probability Measures }

In a predetermined model of uncertainty, the objective is to define, for each event $H$, an interval of real numbers within $[0,1]$ (potentially degenerate). These intervals are intended to represent the degree of confidence, in relation to a probability measure $P$ (or more generally, a finitely additive measure), concerning the occurrence of a specific random phenomenon determined by the causes associated with event $H.$

\smallskip  For this purpose, if a weak complementation $\phi$ is considered on $\mathcal{A}$, the set $H \cup H_{\psi}$ represents the most extensive event that can capture the occurrence of the random phenomenon when the causes are associated with $H$. Moreover, when considering $({\mathcal A},\cup_{\psi_{\mathcal Z}},\cap_{\psi_{\mathcal Z}},\subseteq_{\psi_{\mathcal Z}},{\psi_{\mathcal Z}})$, since $H \cup H_{\psi_{\mathcal Z}}=H\cup_{\psi_{\mathcal Z}} \emptyset$ due to $(12)$, the concept of interval probability measure becomes dependent on both $H$ and $H\cup_{\psi_{\mathcal Z}} \emptyset$. Consequently, it is natural to introduce the following definition, which does not rely on the specific properties of the algebraic structure of ${\mathcal A}$.

\begin{definition}
The map $P_{\psi}$ defined on ${\mathcal A}$ by 
\begin{equation}
P_{\psi}(H)=\Big[P(H),P(H)+P(H_\psi)\Big]
\end{equation}
is the interval probability measure according to $P$ and relative to (the uncertainty determined by) $\psi$. The interval $P_{\psi}(H)$ is said to be the measure (according to $P$) of $H$ relative to ${\psi}$.
\end{definition}
 
\smallskip It is immediately observed that the width $|P_\psi(H)|$ of the interval $P_\psi(H)$ is the probability according to $P$ of the event $H_{\psi}$. Therefore, the interval measure of $H$ coincides with the classical one
according to $P$ if and only if the set of uncertain opportunities of $H$ (i.e., the uncertainty induced by the causes determining $H$ and realizing
the random phenomenon) is $P$-negligible.

\smallskip Moreover, thanks to condition $(b)$ of the weak complementation $\psi$, it follows the decreasing of the indecisiveness of $P_\psi$, that is 
\begin{equation}
H \subseteq K \Longrightarrow |P_\psi(K)| \leq |P_\psi(H)|
\end{equation}

\smallskip\noindent It is important to observe that $P(H) + P(H_\psi) = 1 - P(\psi(H))$. Consequently, the right-end extreme of $P_{\psi}(H)$ aligns with the left-end extreme of $P_{\psi}((\psi(H))^c)$ but not necessarily with $P(\psi(\psi(H)))$ ($\subseteq P(H)$ if $\psi$ is a regular complementation). Before determining the properties of $P_\psi$, a general class of interval probability measures or imprecise probabilities, containing $P_\psi$, is introduced.
\begin{definition}
 An interval
probability measure (or imprecise probability) $Q$ is an application from ${\mathcal{A}%
}$ to the the closed subintervals of $[0,1]$, satisfying

\smallskip\smallskip $(\alpha)$\ \ The left-end extreme $Q_l$ of $Q$ is a probability measure; 

\smallskip $(\beta)$\ \ For every pair of events $H,K$, with $H\subseteq K$, it holds $|Q(K)|\leq |Q(H)|$.

\end{definition}

\bigskip\noindent Condition $(11)$ implies that $P_\psi$ is an interval probability measure in the sense of Definition 6, regardless of whether $\psi$ is a weak complementation. Moreover, there exist interval probability measures $Q$ that do not take the form $P_{\psi}$ for every weak complementation $\psi$, as demonstrated in the following example.

\smallskip\noindent
\begin{example}
Let $Q$ be the interval probability measure defined by $$Q(H)=[P(H), P(H) + rP(H^c)]$$where $P$ is a discrete probability measure and $r \in (0, 1)$ (Obviously, $Q$ verifies conditions $(\alpha)$ and $(\beta)$). The parameter r can be considered as a degree of reliability of uncertain representation, that is a sort of \lq\lq weight of the argument\rq\rq  (Keynes 1921). If $Q$ coincides with $P_{\psi}$ for some weak complementation $\psi$, then $P(H_\psi) = rP(H^c)$ for any event $H$. However, this is not possible when $H^c$ is a singleton because $H_\psi\subseteq H^c$. A similar example holds when $P$ is not a discrete probability measure.
\end{example}

\medskip\noindent An important class of imprecise probabilities can be introduced when  $(rI_{H^c})_{H\in{\mathcal A}}$ of example 8 is substituted by a suitable family $(Y_H)_{H\in{\mathcal A}}$ of random variables. More precisely, 
\begin{definition}
A family $(Y_H)_{H\in{\mathcal A}}$ of random variables, with $Y_H(\Omega)\subseteq [0,1]$ and $Y_H(H)=0$, is said to be not increasing if $Y_K\leq Y_H$ for every pair of events $H,K$, with $H\subseteq K$. 
\end{definition}

\noindent The random variable $Y_H$ can be interpreted as representing the degree of uncertainty associated with the valuation of $H$. When a weak complementation $\psi$ is introduced, it becomes evident that the family $(Y_H)_H$ defined by $Y_H=I_{H_\psi}$ constitutes a non-increasing family of uncertainties. In essence, this implies that the concept of weak complementation can be generalized by incorporating the notion of a non-increasing family of uncertainties.

\medskip\noindent
\begin{definition}
An interval probability measure $Q$ is associated with a not increasing family ${\mathcal Y}=(Y_H)_H$ and a probability measure $P$ if 
\begin{equation}
Q(H) = \big[P(H), P(H) + E^{P}[Y_{H}]\big],
\end{equation}
where $E^{P}[Y_{H}]$ is the mean value of $Y_H$ according to $P$. In this case, $Q$ is denoted by $P_{\mathcal Y}$.
\end{definition}

\begin{remark} Definition 9 implies that the interval measure $Q$ introduced in Example 8 belongs to the class $P_{(rI_{H_\psi})_H}$ and, more precisely, $Q=P_{\mathcal Y}$ when $\mathcal Y=(rI_{H^c})_H$. It is noteworthy that as $r$ increases, the width of $Q$ also increases (i.e., the weight of the argument is heightened, as discussed by Keynes in 1921, p.77). Additionally, the class of interval measures satisfying (12) is stable under convex combinations of $P$, but it may not necessarily be stable under convex combinations involving both $P$ and $(Y_H)_{H\in{\mathcal A}}$.
\end{remark}

\smallskip
\begin{example} Let $M \in \mathbb{N} \cup \{\infty\}$ and $\{\alpha_i\}_{i=0}^M$ be a sequence of mappings defined on a $\sigma$-algebra $\mathcal{A}$ with values in $[0,1]$, such that $\alpha_i(H) \geq \alpha_i(K)$ whenever $H \subseteq K$. Additionally, let $\{K_i\}_{i=0}^M$ be a partition of $\Omega$. In this context, a non-increasing family $\mathcal{Y} = (Y_H)_H$ is constructed by defining
\begin{equation}
Y_H = I_{H^c} \sum_{i=0}^M \alpha_i(H) I_{K_i},
\end{equation}
which represents the uncertainty associated with the valuation of $H$. Notably, if $\alpha_i(H) = I_{\{K_i \subseteq H^c\}}$ and $\mathcal{Z} = \{K_0, \ldots, K_M\}$, this yields $Q_{\mathcal Y}= P_{\psi_{\mathcal{Z}}}$ because $Y_H$ coincide with $I_{H^c\setminus\psi_{\mathcal Z}(H)}$. Furthermore, given a probability measure $P$, the interval probability measure $P_{\mathcal{Y}}$ is defined by
\begin{equation}
P_{\mathcal{Y}}(H) = \big[P(H), P(H) + \sum_{i=0}^M \alpha_i(H)P(H^c \cap K_i)\big].
\end{equation}
\end{example}

\medskip\noindent The interval probability measures $P_{\mathcal{Y}}$ and, more generally, their convex combinations, varying both $P$ and $\mathcal{Y}$, constitute arguably the most natural class of interval probability measures to consider in applications involving Imprecise Probabilities. Nevertheless, let us now delve into the examination of some properties of the interval probability measure $P_{\psi_{\mathcal{Z}}}$, where $\mathcal{Z}$ is a given partition $\{Z_0, \ldots, Z_M\}$ of $\Omega$. Although the measure $P_{\psi_{\mathcal{Z}}}$ may not be overly complicated, its properties are nontrivial and hold a certain element of surprise.

\begin{proposition}
Let $Q$ be an interval probability measure, with left-end extreme $Q_l$. If $Q$ satisfies the following conditions

\smallskip $(i)$ The emptyset and the class of events $H\in{\mathcal{A}}$ such that $%
Q(H)=[a,1]$ and $Q(H^c)=[1-a,1]$, with $a\in[0,1]$, forms a finite
subalgebra ${\mathcal{W}}$ of ${\mathcal{A}}$,

\smallskip $(ii)$ For every atom $Z\in{\mathcal{W}}$ and for any non-empty
event $H\subseteq Z$, it holds 
\[
Q(H)=[Q_l(H),1-Q_l(Z\cap H^c)], 
\]

\smallskip $(iii)$ For every pair $Z,W$ in ${\mathcal{W}}$, with $Z$ atom
disjoint from $W$, and for any pair $H,K$ of non-empty events respectively
contained in $Z,W$, it results 
\[
|Q(H\cup K)|+Q_l(Z)=|Q(K)|,
\]

 \noindent then  $Q=P_{\psi_{\mathcal{Z}}}$, where $P=Q_l$ and $\mathcal{Z}$ is
determined by the atoms of $\mathcal{W}.$ In particular, if ${\mathcal{W}}$ does not contain any proper subset of $\Omega$, then  $Q(H)=[Q_l(H), Q_l(H)]$ for every non-empty event in ${\mathcal{A}}$.  

\medskip Conversely, if $Q=P_{\psi_{\mathcal{Z}}}$ for some finite partition $\mathcal{Z}$, then conditions $(i)-(iii)$ are fulfilled, where ${\mathcal W}$ is the algebra generated by ${\mathcal Z}$. 
\end{proposition}

\medskip\noindent Conversely, isotonic properties of the type $H\subset
K\Longrightarrow P_{\psi_{\mathcal Z}}(H)\subseteq P_{\psi_{\mathcal Z}}(K)$ or
endpoints of $P_{\psi_{\mathcal Z}}(H)$ are rispectively smaller than $P_{\psi_{\mathcal Z}}(K)$ are not true.
Moreover, when $P$ is a probability
measure, for every increasing sequence $(H_n)_n$ of events with $\cup_nH_n=H$, it holds 
\begin{equation}
\lim_n P((H_n)_{\psi_{\mathcal Z}})=P(H_{\psi_{\mathcal{Z}}}).
\end{equation} Then the endpoints of $P_{\psi_{\mathcal Z}}((H_n)_{\psi_{\mathcal Z}})$ converge to that of $P_{\psi_{\mathcal Z}}(H_{\psi_{\mathcal{Z}}})$;
the same holds if $(H_n)_n$ is a decreasing sequence of events and $P_{\psi_{\mathcal Z}}(H)$ is a degenerate interval. 

\smallskip It is worth noting that the
same statement does not hold when $P_{\psi_{\mathcal{Z}}}(H)$ is not a degenerate interval.
Specifically, it suffices to consider a decreasing sequence with an empty
intersection and all $H_n$ intersecting each $Z_k$. Notice the different
behavior of the monotonicity property between decreasing sequences of events
and increasing sequences of events.

\begin{proposition}
Let $Q$ be an interval probability measure and $\Omega=\{\omega_i\}_{i\in J}$ be a countable space, with $Q_l(\{\omega_i\})\neq 0$ for every $i\in J$.
Then $Q=P_{\mathcal Y}$ for a suitable not increasing family ${\mathcal Y}$ and $P=Q_l$ iff for every finite subset $I$ of $J$ it holds $|Q({\Omega\setminus\{\omega_i\}_{i\in I}})|\geq \sum_{i\in I}|Q(\Omega\setminus\{\omega_i\})| $.
\end{proposition}

\smallskip\noindent
\begin{remark} In place of the probability
measure $P$, a more generalized definition of interval measure $\nu_{\psi_{\mathcal{Z}}}$
can be derived by considering a capacity $\nu$, namely an application $\nu$
defined on ${\mathcal{A}}$ such that $%
\nu(\emptyset)=0, \nu(\Omega)=1$ and $\nu(H)\leq\nu(K)$ when $H\subseteq K$. More
precisely, $\nu_{\psi_{\mathcal{Z}}}$ could be defined as follows 
$$
\nu_{\psi_{\mathcal{Z}}}(H)=[\nu(H),\nu(H\cup H_{\psi_{\mathcal{Z}}})]. 
$$
It is evident that several properties, such as $(ii)$ and $(13)$, are no
longer satisfied. Another potentially interesting generalization arises
when $Z_n \notin {\mathcal{A}}$, and $P_{\psi_{\mathcal Z}}$ is defined as
follows: 
\[
P_{\psi_{\mathcal Z}}(H)=[P(H),P(H)+\sup_{B\subseteq H_{\psi_{\mathcal{Z}}},B\in{%
\mathcal{A}}}P(B)]. 
\]
These two generalizations will be explored in an upcoming paper.
\end{remark}

\medskip In the following example, we explore a finite collection of interval probability measures and consequently establish a representation of the overall uncertainty associated with the collection.

\begin{example}
Consider $h$ interval probability measures $(Q_j)_{j=1,\ldots,h}$, where $Q_j$ denotes $(P_{j})_{{\mathcal Y}_j}$ and ${\mathcal Y}_j=(Y_{j,H})_H$ is not-increasing family of random variables. The imprecise probability of $H$ relative to the
collection $(Q_j(H))_j$ is defined by
\[
\widetilde Q(H)=\Big[\inf_j P_j(H), \sup_j E^{P_j}[I_H+Y_{j,H}]\Big]%
. 
\]
However, $\widetilde Q$ is not an interval probability measure because it does not satisfy property $(\beta)$. Moreover, if 
$Y_{j,H}$ is defined as in $(16)$, where $(K_{i,j})_i$ is a partition of $\Omega$ with $P_j(K_i)\neq 0$
and $\alpha_{i,j}(H)\leq P_j(K_i)$ for any $i$, the global
uncertainty associated with $(Q_j(H))_j$ can be defined as 
\[
Y_H=I_{H^c}\sup_j\sum_{i=1}^\infty \frac{\alpha_{i,j}(H)}{P_j(K_{i,j})%
}I_{K_{i,j}}. 
\]
In such a case, an other imprecise probability measure of $H$ can be defined
as 
\[
\widehat Q(H)=\Big[\inf_j P_j^\prime (H),
\sup_j(P_j^\prime(H)+E_{P^\prime_j}[Y_H])\Big]. 
\]
It is immediately observed that $\widetilde Q(H)\subseteq\widehat Q(H)$ for
any event $H$. A measure of the coherence of uncertainties relative to $(Q_j(H))_j$ could be obtained by considering 
\[
\delta(H)=\sup_j\big(1-{\frac{{|Q_j(H)|}}{{E_{P_j}[Y_H]}}}\big)
\]
For small values of $\delta(H)$ the coherence of uncertainties is high. In the particular case where $Y_{j,H}=I_{{H_{\psi_{{\mathcal Z}_j}}}}$, it follows $Y_H=I_{\cup_j{H_{\psi_{{\mathcal Z}_j}}}}$ and $$\delta(H)= \sup_j\big(1-\frac{{P_j(H_{\psi_{{\mathcal Z}_j}})}}{{P_j}(\cup_{i=1}^h H_{\psi_{{\mathcal Z}_i}})}\big) $$ High coherence in uncertainties is achieved when the partitions ${\mathcal Z}_j$ are very similar. This
approach can be applied when aggregation of assessments and judgments is
required.
\end{example}

\medskip\noindent Concluding this section, we finally provide a simple example that constructs a partition ${\mathcal Z}$ and an interval probability measure $P_{\psi_{\mathcal Z}}.$

\begin{example}
Consider a slightly modified version of the "Treatment of Uncertainty" used
in the IPCC's Assessment Reports 5 (AR5) and IPCC's Sixth Report (2021). Let
the set of ambiguous or vague belief be as follows: \textit{Virtually
Certain; Very Likely, Likely, About as Likely as Not, Very Unlikely and
Exceptionally Unlikely}.

\smallskip It is reasonable to assume that an agent possesses a subjective 
\textit{belief function} (and consequently, a probability measure) for an increasing class of events, reflecting her degrees of belief in those events. As a
result, an ordinal scale of uncertain beliefs can be associated with a scale
of probability intervals, such as: [98\%,100\%]; [90\%,98\%]; [66\%,90\%];
[33\%,66\%]; [10\%,33\%]; [2\%,10\%]; [0\%,2\%].
\end{example}

In our probabilistic model, this can be achieved by taking $M=6$ and
assuming a partition ${\mathcal{Z}}$ of $\Omega$ into seven events $Z_0,\ldots,Z_6$ along with a probability measure $P$ such that 
\[
P(Z_6)=P(Z_0)=2\%,\quad P(Z_5)=P(Z_1)=8\%, \quad P(Z_4)=24\%, 
\]
\[
P(Z_3)=33\%,\quad P(Z_2)=23\%; 
\]
considering events $H_n=\cup_{i=0}^{n-1}Z_j\cup_{j=n+1}^6 W_j,$ where $W_j$
is a non-empty negligible event (or $P(W_j)$ is much smaller than $0.08\%)$
contained in $Z_j$, it holds 
\[
P_{\psi_{\mathcal{Z}}}(H_n)=[\sum_{i=0}^{n-1} P(Z_i),\sum_{i=0}^{n}P(Z_i)] 
\]
and, in particular, the following interval probability measures are obtained 
\[
P_{\psi_{\mathcal{Z}}}(H_6)= [98\%,100\%], P_{\psi_{\mathcal{Z}}}=[90\%,98\%], P_{\psi_{\mathcal{Z}}}(H_4)=[66\%,90\%] 
\]
\[
P_{\psi_{\mathcal{Z}}}(H_3)= [33\%,66\%], P_{\psi_{\mathcal{Z}}}(H_2)=[10\%,33\%], P_{\psi_{\mathcal{Z}}}(H_1)=[2\%,10\%] 
\]
and $P_{\psi_{\mathcal{Z}}}(H_0)=[0\%,2\%].$

\medskip

\section{Conditional Interval Probability Measures}

Let $Q$ be an interval probability measure associated with a not increasing family ${\mathcal Y}=(Y_H)_H$. According to Definition 10, $Q$ coincides with $P_{\mathcal Y}$ for a suitable probability measure $P$. In the following, let $H$ be an event not negligible under $Q$—that is, $Q(H)\neq [0,0]$. An interval probability measure can be defined to express the degree of confidence in realizing the random phenomenon conditioned on the causes that determined $H$. For this purpose, for each event $A$ in ${\mathcal{A}}$, we define the \textit{probability of $A$ conditioned on $H$ relative to $Q$} as
\begin{equation}
Q(A|H)=\big[{\frac{{P(A\cap H)}}{{P(H)+E^{P}[Y_H]}}},{\frac{{E^{P}[(I_A+Y_A)(I_H+Y_H)]}}{{P(H)+E^{P}[Y_H]}}}\big].
\end{equation}
If $Q=P_{\psi_{\mathcal{Z}}}$ where ${\mathcal Z}$ is a finite partition, it leads to
\begin{equation}
P_{\psi_{\mathcal{Z}}}(A|H)=\Big[P\big(A\, | (\psi_{\mathcal{Z}}(H))^c\big),P\big((\psi_{\mathcal{Z}}(A))^c\, |(\psi_{\mathcal{Z}}(H))^c\big)\Big] 
\end{equation}
It is noteworthy that $P_{\psi_{\mathcal Z}}$ possesses left and right extremes that adhere to the duality rule concerning weak complementation
\begin{equation}
P\big((\psi_{\mathcal{Z}}(A))^c\, |(\psi_{\mathcal{Z}}(H))^c\big)=1-P\big(\psi_{\mathcal{Z}}(A)\, |(\psi_{\mathcal{Z}}(H))^c\big).
\end{equation} Moreover, it holds
$$P_{\psi_{\mathcal{Z}}}(\emptyset|H)=[0,1] \quad P_{\psi_{\mathcal{Z}}}(\Omega|H)=P_{\psi_{\mathcal{Z}}}((\psi_{\mathcal{Z}}(H))^c\, |H)=[1,1] 
$$
while $P_{\psi_{\mathcal{Z}}}(H\, |H)=[P(H)/(1-P(\psi_{\mathcal{Z}}(H))),1]$.

\smallskip\noindent

\begin{proposition}
Let $H$ be an event not  $P_{\psi_{\mathcal{Z}}}$-negligible and $H\subseteq Z_n$
for some $n$. If $A$ is an event contained in some $Z_m$, for $m\neq n$,
then 
\begin{equation}
P_{\psi_{\mathcal{Z}}}(A|H) = \left[ \frac{P(A)}{P(H) + 1 - P(Z_n)}, \frac{P(H) + 1
- P(Z_n) - P(Z_m\cap A^c)}{P(H) + 1 - P(Z_n)} \right]
\end{equation}
while, for $n=m$, it results 
\begin{equation}
P_{\psi_{\mathcal{Z}}}(A|H) = \left[ \frac{P(A\cap H)}{P(H) + 1 - P(Z_n)}, \frac{%
P(A\cap H) + 1 - P(Z_n) }{P(H) + 1 - P(Z_n)} \right]. 
\end{equation}
In particular, $P_{\psi_{\mathcal{Z}}}(Z_m|Z_n) = \big[P(Z_m), 1\big] = P_{\psi_{\mathcal{Z}}}(Z_m),$ for any $m,n$.
\end{proposition}

\medskip\noindent Additionally, if $\psi_{\mathcal{Z}}(H)=H^c$, the probability of $A$ conditioned on $H$ relative to $P_{\psi_{\mathcal{Z}}}$ is the interval $[P(A|H),P(A|H)+P(({\psi_{\mathcal{Z}}}(A))^c)|H)]$. In this case, $%
P_{\psi_{\mathcal{Z}}}(A|H) $ is degenerate if and only if $P(({\psi_{\mathcal{Z}}}(A))^c\cap
H)=0$.

\medskip\noindent \textbf{Example 17 (continued)} For any pair $m,n$ with $%
m<n\leq 6 $, it holds 
\begin{equation}
P_{\psi_{\mathcal{Z}}}(H_m|H_n)=\Big[{\frac{{\sum_{i=0}^{m-1} P(Z_i)}}{{%
\sum_{i=0}^{n} P(Z_i)}}},{\frac{{\sum_{i=0}^{m} P(Z_i)}}{{\sum_{i=0}^{n}
P(Z_i)}}}\Big]
\end{equation}
Thus, the 'knowledge' of the causes determining $H_n$ increases the
uncertainty induced by the causes determining $H_m$. Consequently, $H_m$ is
not considered independent of $H_n$, as we will explore in the next section.

\medskip
\begin{remark}
Let $\nu$ be a super additive capacity  i.e. for any event $H$, $\nu$ satisfies $\nu(H)+\nu(H^c)\leq 1$.  If $\nu(H^c)\neq 1$, Dempster and Shafer define  $$\nu(A|H)={{\nu((A\cap H)\cup H^c)-\nu(H^c)}
\over{1-\nu(H^c)}}$$ for any event $A$. 
This definition does not equate to standard Bayesian updating if $\nu$ is not additive. 

However, if $H^c$ is replaced with the weak complementation $\psi_{\mathcal Z}(H)$, and thus $H$ replaced with $(\psi_{\mathcal Z}(H))^c$, interpreting the plausibility $1 - \nu(\psi_{\mathcal Z}(H))$ as the measure of $(\psi_{\mathcal Z}(H))^c$ according to a probability $P$, we can see that $${{\nu((A\cap (\psi_{\mathcal Z}(H))^c)\cup\psi_{\mathcal Z}(H))-\nu(\psi_{\mathcal Z}(H))}
\over{1-\nu(\psi_{\mathcal Z}(H))}}={{\nu((A\cup \psi_{\mathcal Z}(H))-\nu(\psi_{\mathcal Z}(H))}
\over{1-\nu(\psi_{\mathcal Z}(H))}}$$ can be understood as the ratio of measures of $A \cap (\psi_{\mathcal Z}(H))^c$ and $(\psi_{\mathcal Z}(H))^c$ according to $P$. In other words, the left endpoint of $P_{\psi_{\mathcal F}}(A|H)$. Similarly, the interpretation of the right endpoint of $P_{\psi_{\mathcal F}}(A|H)$ holds true, thanks to duality equation $(21)$. Therefore, for a suitable $P$, the probability $P_{\psi_{\mathcal F}}(A|H)$ can represent the interval where the extremes are the two conservative belief and plausibility degrees conditional to $H$ by Dempster and Shafer.
\end{remark}

\medskip Are there other 'reasonable' ways to define the concept of conditional probability with respect to $H$ according to a general interval probability measure $Q$? It is quite likely. Perhaps one might consider situations where uncertainty is negligible with respect to $Q$ (i.e., $|Q(H)|=0$) and define
\[
Q(A|H)=\big[Q_l(A|H),{\frac{{Q_r(A\cap H)}}{{Q_l(H)}%
}}\big]. 
\]
However, it seems to us that the concept of conditional interval probability is coherent only when considering interval probability measures $Q=P_{\mathcal{Y}}$. For now, let's continue exploring the consequences of the notion of conditional probability with respect to $P_{\psi_{\mathcal{Z}}}$.

\medskip

\section{Independence of an event $A$ conditioned to $H$}

With the same assumptions and notations as in the previous section, we will say that an event $A$ is
independent of $H$ according to $Q=P_{\mathcal{Y}}$ if and only if 
\begin{equation}
Q(A|H)=Q(A)
\end{equation}
Following relation $(20)$, it immediately follows that there is
independence of $A$ from $H$ according to $P_{\psi_{\mathcal{Z}}}$ if and only if $A $ and $\psi_{\mathcal{Z}}(A)$ are  independent from $(\psi_{\mathcal{Z}}(H))^c$ with respect to $P$.  

\smallskip Note that for an interval probability measure,
the notion of an event $A$ being independent of an event $H$ does not
necessarily coincide with the independence of $H$ from $A$. To illustrate this 'asymmetry' in the concept of independence, consider $Q=P_{\psi_{\mathcal Z}}$ where ${\mathcal Z}$ is the partition $Z_0,Z_1$. Moreover, let $H$ be a not $P$-negligible event which intersects $Z_0$ and $Z_1$, but $P(H\cap Z_1)=0$. If $P(Z_1)>0$ then $H$
is independent of $Z_0$ (according to $P_{\psi_{\mathcal Z}}$) because 
\[
P_{\psi_{\mathcal Z}}(H|Z_0)=P_{\psi_{\mathcal Z}}(H)=[P(H),P(H)] 
\]
but $Z_0$ is not independent from $H$ (according to $P_{\psi_{\mathcal Z}}$) since $P_{\psi_{\mathcal Z}}(Z_0)=[P(Z_0),1]$ and $P_{\psi_{\mathcal Z}}(Z_0|H)=[1,1]$.

\smallskip Moreover, owing to Proposition 18, it turns out that the events $%
Z_0$,$Z_1$ are independent in the sense that both $Z_0$ is 
independent of $Z_1$ and vice versa. Similarly, perhaps
paradoxically, $Z_0$ and $Z_0$ as well as $Z_1$ with $Z_1^c$ are independent.

\smallskip In simpler terms, altering an event on a set $B$ that is negligible under the probability measure $P$, but has a non-degenerate interval probability measure $P_{\psi_{\mathcal{Z}}}(B)$, may transition 
from $P_{\psi_{\mathcal{Z}}}$-independence to non-$P_{\psi_{\mathcal{Z}}}$-independence.
Similar considerations apply when considering a general $P_{\mathcal{Y}}$.

\medskip\noindent \textbf{Example 17 (continued)} The previous apparent
paradox it becomes evident when considering $H_m, H_n$ and $%
\cup_{i=1}^{m-1}Z_i$, $\cup_{i=1}^{n-1}Z_i$ with $1\leq m<n\leq 6$. In fact 
\[
P(H_m\neq \cup_{i=1}^{m-1}Z_i)=0,\qquad P(H_n\neq \cup_{i=1}^{n-1}Z_i)=0 
\]
and $H_m$ is not independent of $H_n$ (according to $P_{\psi_{\mathcal Z}}$) thanks to $(23)$. However, $\cup_{i=1}^{m-1}Z_i$ and $\cup_{i=1}^{n-1}Z_i$ are independent because 
\[
P_{\psi_{\mathcal Z}}(\cup_{i=1}^{m-1}Z_i|\cup_{i=1}^{n-1}Z_i)=\Big[%
\sum_{i=1}^{m-1}P(Z_i),1\Big]=P_{\psi_{\mathcal Z}}(\cup_{i=1}^{m-1}Z_i) 
\]

\smallskip

\section{Interval Distribution of a Random Variable}

Now, let's consider a real-valued random variable $X$. The \textit{%
distribution function according to $P_{\psi_{\mathcal{Z}}}$ of $X$} is defined as
the function $F_{\mathcal{Z}}$ from the real line to the set of closed
subintervals of $[0,1]$ given by 
\begin{equation}
F_{\mathcal{Z}}(t)=P_{\psi_{\mathcal{Z}}}(X\leq t)
\end{equation}
In general, if an interval probability measure $Q$ is considered the
distribution function is defined by $F_{Q}(t)=Q(X\leq t).$

\begin{example}
A special case arises when $X(Z_i)=t_i$, $Q=P_{\psi_{\mathcal{Z}}}$ and thus $X$
takes at most countably many values. In this case, it holds that 
\[
F_{\mathcal{Z}}(t)=\Big[P(X\leq t),1\Big]
\]
where $P(X\leq t)=\sum_{i:t_i\leq t} P(Z_i)$. More generally, if $t_i$ is non-decreasing and $
Y(Z_i)=]t_{i-1},t_i]$, then the distribution function $G_{\mathcal{Z}}$,
according to $P_{\psi_{\mathcal{Z}}}$, for the continuous random variable $Y$ is
given by 
\[
G_{\mathcal{Z}}(t)=\Big[P(Y\leq t), P(Y\leq t) + \sum_{i:t_i>t}P(Z_i)\Big], 
\]
meaning that the endpoints of the interval confidence are the distribution
function $G$ of $Y$ according to $P$ and the sum of $G$ and the decreasing
function 
$$
t\mapsto \sum_{i:t_i>t}P(Z_i). 
$$
In particular, $
G_{\mathcal{Z}}(t)=\big[G(t), 1-P(t<Y\leq t_{*})\big]=\big[G(t), 1-G(t_*) +
G(t)\big]$, where $t_*$ is the unique $t_i$ such that $t\in]t_{i-1},t_i]$. Note that the
function representing the right endpoint of $G_{\mathcal{Z}}$ is neither
increasing nor decreasing. 

\smallskip Moreover, the widths of the distribution function 
$F_{\mathcal{Z}}$ and $G_{\mathcal{Z}}$ depend significantly on the images
of $Z_i$. Observe that if $X\geq Y$, it results $F_{\mathcal{Z}%
}(t)=G_{\mathcal{Z}}(t)$ when $t=t_i$ for any $i$, while $F(t)\leq G(t)$ and 
$|F(t)|\geq |G(t)|$ if $t\neq t_i$.
\end{example}

\medskip \noindent The previous example suggests a notion of stochastic dominance with respect to an interval probability measure $Q$.

\begin{definition}
$X$ stochastically dominates $Y$ with respect to $Q$ if 
\begin{equation}
Q_{l}(X\leq t)\leq Q_{l}(Y\leq t) \ \ \ {and}\ \ \ |Q(X\leq t)|\geq
|Q(Y\leq t)|\quad \ \ \ \ \forall t\in {\mathbb{R}}, 
\end{equation}
where $Q_{l}$ represents the left-end extreme of $Q$.
\end{definition}

 \noindent A similar definition can be extended when $X$ and $Y$ are random vectors. This concept of interval stochastic dominance enables the comparison of random variables in terms of uncertainty and provides a straightforward guideline for making choices among various alternatives.

\smallskip Finally, for a random variable $X$ with values in a measurable
space $(E,{\mathcal{E}})$, the {\sl law of $X$ according to $Q$} is the
interval measure $\mu=X(Q)$ defined on ${\mathcal{E}}$ by 
$$
\mu(H) = Q\big(X^{-1}(H)\big). 
$$
Note that 
\[
P_{\psi_\mathcal{Z}}\big(X^{-1}(H)\big) = [\mu(H), \mu(H) + \sum_{n:X^{-1}(H)\cap
Z_n= \emptyset} P(X^{-1}(H)\cap Z_n)], 
\]
where $\mu$ denotes the usual law of $X$ according to $P$. If $%
Z_n=X^{-1}(K_n)$ then 
\[
P_{\psi_\mathcal{Z}}\big(X^{-1}(H)\big) = [\mu(H), \mu(H) + \sum_{n:H\cap K_n
=\emptyset} \mu(H\cap K_n)], 
\]

\medskip\noindent
\section{A possible solution of controversy between \newline
Key\-nes and Ramsey.}

In \textit{Chapter 2} of \textit{Truth and Probability}. Ramsey sets that
the degree of the probability relation is the same as the the degree of
belief and that both of them are expressed by numbers that are the same
(1926, 160). Ramsey assumes that there is the one-to-one relation between
degree of belief and probability. As a consequence, the one-to-one
relation between belief and probability is an order bijective function or
isomorphism, that is the function clearly preserves the relations of greater
and less among elements of the sets.

Inexplicably, for a century, up to now the seminal problem put in
evidence by Ramsey was never considered. Nevertheless, it is the formal core
of the Ramsey's doubt about the Keynes' notion of non-numerical probability.

Following the seminal paper by Birkhoff and von Neumann (1936) on quantum mechanics, a lattice is defined, which, unless considering isomorphisms, corresponds to an algebra of sets. In such an algebra, if a (finite) partition ${\mathcal Z}$ of events is given, we can identify for any event $H$ the set of vague belief with the interval $[P(H),P((\psi(H))^c)]$, representing the range of approximate probabilities.

Moreover, it can be demonstrated that the formal structure of Ramsey's theory and Keynes's theory is common, but they differ in terms of set complementation. For Keynes, the complementation of $H$ can be expressed as $
\psi_{\mathcal Z}(H)$, and Keynesian uncertainty is represented by 
$H^c\cap (\psi_{\mathcal Z}(H))^c$, while for Ramsey, the complementation of $H$ is always the usual set $H^c$. It follows that the isomorphism is unique only if ${\mathcal Z}$ is determined or fixed.

In Section 3, we introduce the notion of interval probability measure, revealing that there are infinite interval probability measures of the type $P_{\psi_{\mathcal Z}}$ in Keynes' pseudo-Boolean algebra that can represent uncertain beliefs. Importantly, all these measures $P_{\psi_{\mathcal Z}}$ may describe the evaluation of the uncertainty of an event, even when Ramsey's probability measure $P$ is already given.

\smallskip Consider the famous example of the umbrella: {\sl  Is our expectation of rain,
when we start out for a walk, always more likely than not, or less likely
than not, or as lkely as not? I am prepared to argue that on same occasions
none of these alternatives hold, and that it will b an arbitrary mater to
decide for or against the umbrella. If the barometer is high, but the clouds
are black, it is not always rational that one should prevail over the other
in our mind, or even that we should belace them.} (Keynes 1921, 28).

In our framework of imprecise probabilities or interval probability measures, we can represent the problem of the umbrella as a situation where there are four relevant causes for raining:

\smallskip\noindent
$i)$  barometer low and none black clouds in the sky, denoted as eventuality $\omega_{01}$;

\noindent $ii)$ barometer high, many black clouds in the sky, denoted as eventuality $\omega_{10}$;

\noindent $iii)$ barometer low, many black clouds in the sky, denoted as eventuality $\omega_{00}$;

\noindent $iv)$ barometer high, none black clouds in the sky, denoted as eventuality $\omega_{11}$.

\medskip\noindent Following Example 5, since the negation of a given cause does not necessarily imply the non-occurrence of the represented event, we establish the partition  ${\mathcal Z}=Z_0,Z_1$ with $Z_0=\{\omega_{0,1},\omega_{1,0}\}$, $Z_1=\{\omega_{0,0},\omega_{1,1}\}$ to describe the weak complementation and, consequently, uncertain events.
In particular, the Law of
the Exclude Middle does not held because the (weak) complementary event of $\{\omega_{10}\}$ is $\{\omega_{01}\}$
and $Z_1$ represents the set of indecisive eventualities of $\{\omega_{10}\}$.

\smallskip In a nutshell, when facing discordant situations like $\omega_{10}$ and $\omega_{01}$, it is not possible to
infer anything about whether these are due to a distorted or mis-revealed occurrence of agreeing causes in $Z_1$. Similarly, when a decision-maker encounters $\omega_{11}$ and $\omega_{00}$, which are agreeing situations, it does not allow them to deduce whether these are mis-reliable occurrences of discordant situations ($\omega_{10}$ and $\omega_{01}$) or
not. In our framework, if $P$ is a probability on the subset of $\Omega=\{\omega_{ij}\}_{i,j=0,1}$ and $H=\{\omega_{10}\}$, 
the interval probability of $H$ is given by
$$P_{\psi_{\mathcal Z}}(H)=[P(H),P(H)+P(Z_{1})]=[P(H),1-
P(\{\omega_{01}\})]$$ while
$$P_{\psi_{\mathcal Z}}(\{\omega_{11}\}|H)=[{{P(\{\omega_{11}\})}\over{1-P(\{\omega_{01}\})}},{{P(\{\omega_{11}\})+P(H)}\over{1-P(\{\omega_{01}\})}}]$$ and $$P_{\psi_{\mathcal Z}}(\{\omega_{00}\}|H)=[{{P(\{\omega_{00}\})}\over{1-P(\{\omega_{01}\})}},{{P(\{\omega_{00}\})+P(H)}\over{1-P(\{\omega_{01}\})}}].$$
Moreover $$P_{\psi_{\mathcal Z}}(\{\omega_{00}\}\cup\{\omega_{11}\}|H)=[{{P(\{\omega_{00}\})+P(\{\omega_{11}\})}\over{1-P(\{\omega_{01}\})}},1]\supset P_{\psi_{\mathcal Z}}(\{\omega_{00}\}\cup\{\omega_{11}\}) $$
Observe, $P_{\psi_{\mathcal Z}}(\{\omega_{00}\}\cup (\{\omega_{11}\}|H)\neq P_{\psi_{\mathcal Z}}(\{\omega_{00}\}|H)\cup P_{\psi_{\mathcal Z}}(\{\omega_{11}\}|H)$.

\medskip If $P(Z_0)$ is a small value, depending on
relaibility of barometer and glance, the widht of $P_{\psi_{\mathcal Z}}(\{\omega_{00}\})=1-P(Z_0)$ (evaluation of uncertainty) becomes very large. Thus everything may happen, and
decision depends on attitude toward uncertainty.

It is essential to emphasize that when $P(Z_0)$ is small, even if $P(\{\omega_{10}\})$ (or $P(\{\omega_{01}\})$) is greater than $P(\{\omega_{01}\})$ (or $P(\{\omega_{10}\})$) and the average point of $P_{\psi_{\mathcal Z}}(\{\omega_{10}\})$ is on the left side (resp. right side) with respect to $\frac{1}{2}$,
the uncertainty is too substantial to conclusively determine whether the umbrella is the best choice or not. 

\smallskip Similar outcomes can be derived by considering an alternative model of uncertainty, which is obtained through the partition  ${\mathcal Z}^\prime=\{\{\omega_{10}\},\{\omega_{01}\}, Z_1\}$. Generally, when a decision-maker is tasked with deciding for or against using an umbrella, they need to delineate a model for uncertainty. In this case, the model is obtained through a partition in the set of possible events (or a not increasing family of random variables) that allows the introduction of the notion of weak complementation. Once a weak complement is defined, a probability measure can be elicited, and the corresponding interval probability measures can be determined.

\smallskip Our concept of imprecise probability leads to a clear conclusion: since the notion of weak (relative) complementation is not unique, there exists a (finite or infinite) family (or ensemble, following von Neumann) of interval probability measures that can model and assess sets of eventualities for a random phenomenon, even when the left extreme of these measures is specified.

\medskip Consequently, the framework proposed by Ramsey is merely a specific case of Keynes' framework, as the Law of the Excluded Middle is derived by assuming that weak (relative) complementation aligns with set complementation.

In summary, the isomorphism between belief and probability in Keynes' framework cannot be demanded, as it remains undetermined when the notion of weak (relative) complementation is not predefined.

\smallskip
\section{Conclusion}

\smallskip 
In this paper, we have introduced general classes of interval probability measures to represent uncertainty. Several interval probability measures are defined in relation to regular weak-complementations on an algebra of events. The concept of regular weak-complementation, involving straightforward manipulations of the order relation ($\subseteq$) and standard union ($\cup$), allows for the definition of a complemented but non-distributive lattice. These complemented and non-distributive lattices bear a resemblance to the non-distributive and orthocomplemented lattice introduced by Birckoff and von Neumann to represent Quantum Mechanics, exhibiting a coherent definition of interval probability measures.

Remarkably, the concept of interval probability measures reveals that the convex combination of two measures remains an interval probability measure.  However, if these two measures are generated by weak-complementations, the resulting measure may not establish a weak complementation. This holds true as long as the two weak-complementations do not generate the same elementary events.

Our framework of interval measures enables the introduction of a notion of stochastic dominance between alternatives, making it suitable for uncertain comparison in the decision-making process. Moreover, it facilitates the definition of conditional interval measures through a straightforward updating rule, illustrating how new information modifies priors. The updating rule results in a new interval probability measure that adheres to the duality rule relative to weak complementation, establishing a direct connection with Dempster-Shafer's rule. Naturally, the independence between random events is derived, and its asymmetry is established.

Ultimately, the concept of interval probability measures associated with weak complementation allows us to provide a comprehensive and clear answer to Ramsey's question regarding the relationship between belief and probability. Specifically, the representation of uncertainty through interval probability measures demonstrates that if a partition is defined to introduce weak complementation, there exists a non-unique interval probability measure with the same left-end extreme for every (additive) probability considered.

However, when limiting the consideration to interval probability measures generated solely by fixed probabilities and weak complementations, a unique interval measure is precisely determined. This precise determination of imprecise probability aligns with Ramsey's age-old inquiry.

\section{Appendix}

\medskip
\subsection{ Proof Proposition 4}

Condition $(a)$ of Definition 1 is evidently satisfied because $$\psi_{\mathcal Z}(H) = \bigcup_{n : H \cap Z_n \neq \emptyset} H^c \cap Z_n\subseteq H^c$$ while condition $(b)$ is verified since $$H\subseteq K\Longrightarrow \{n:K\cap Z_n=\emptyset\}\subseteq \{n:H\cap Z_n=\emptyset\}\Longrightarrow \kern-1.5mm\bigcup_{n : K \cap Z_n=\emptyset}  Z_n\subseteq\kern-1.5mm \bigcup_{n : H \cap Z_n=\emptyset}  Z_n$$ Moreover, $\psi_{\mathcal Z}$ is a regular weak-complementation, since $$\{n : \psi_{\mathcal Z}(H) \cap Z_n \neq \emptyset\}=\{n : H \cap Z_n \neq \emptyset, H^c\cap Z_n\neq\emptyset\},$$ and it holds $$H^c\cap\psi_{\mathcal Z}(\psi_{\mathcal Z}(H))=H^c\cap[\psi_{\mathcal Z}(H)]^c\cap \bigcup_{n : \psi_{\mathcal Z}(H) \cap Z_n \neq \emptyset}Z_n=\emptyset $$ because $H_{\psi_{\mathcal Z}}\cap \bigcup_{n : H \cap Z_n \neq \emptyset, H^c\cap Z_n\neq\emptyset}Z_n=\emptyset.$ Relation $(5)$, $(6)$, and $(7)$ follow immediately from $(2)$ and De Morgan's rules of standard complementation.
Since $\{n : H \cap Z_n \neq \emptyset\}\cup\{n : K \cap Z_n \neq \emptyset\}$ is the union $$I_1\cup I_2\cup \{n: H\cap Z_n \neq
\emptyset,K\cap Z_n \neq \emptyset\}$$
with $I_1=\{n: H\cap Z_n \neq \emptyset=K\cap Z_n\}$ and $I_2=\{n: K\cap Z_n \neq \emptyset=H\cap Z_n\},$
the thesis follows from
\begin{eqnarray*}
\psi_{\mathcal Z}(H)\cup \psi_{\mathcal Z}(K) &=&\bigcup_{n : H \cap Z_n \neq \emptyset} H^c \cap Z_n\ \cup \bigcup_{n : K\cap Z_n \neq \emptyset} K^c \cap Z_n\\
&=&
J\cup \bigcup_{n: H\cap Z_n \neq
\emptyset\, \wedge\, K\cap Z_n \neq \emptyset}\, (H^c\cup K^c)\cap Z_n
\\
\end{eqnarray*}
where $J=\cup_{n\in I_1} H^c\cap Z_n \, \cup_
{n\in I_2} K^c\cap Z_n.$ 

\medskip
\subsection{ Proof of (11)}

To prove $(11)$ it suffices to show 
$$H\cup_{\psi_{\mathcal Z}}\psi_{\mathcal Z}(H)=\Omega,$$ since $\psi_{\mathcal Z}(\psi_{\mathcal Z}(H))\subseteq H$ and
$$H\cap_{\psi_{\mathcal Z}}\psi_{\mathcal Z}(H)=H\cap\psi_{\mathcal Z}(H)=\emptyset.$$ Then relation $(11)$ follows from  $$H\cup_{\psi_{\mathcal Z}}\psi_{\mathcal Z}(H)=[\psi_{\mathcal Z}(H\cup \psi_{\mathcal Z}(H))]^c=[\psi_{\mathcal Z}(H^c_{\psi_{\mathcal Z}})]^c=\Omega,$$ because $H^c_{\psi_{\mathcal Z}}=\cup_{\{n:H\cap Z_n\neq\emptyset\}}Z_n$ and $\psi_{\mathcal Z}(\cup_{\{n:H\cap Z_n\neq\emptyset\}}Z_n)=\emptyset$.
\bigskip
\subsection{ Proof Proposition 13}
Let ${\mathcal Z}=\{W_0,\ldots,W_M\}$ be the partition generated by atoms of ${\mathcal W}$. If $M=0$ relation $(ii)$ implies $Q(H)=[Q_l(H),Q_l(H)]$ for any non empty-event and thesis is satisfied. Now, assume $M\geq 1$. To establish $Q=P_{\psi_{\mathcal Z}}$, it is sufficient to demonstrate, for every $m=0,\ldots,M$, that:

\smallskip {\sl For every event $H\subseteq \cup_{k=0}^mZ_k$, with $H\cap Z_k\neq \emptyset$ ($k=0,\ldots, m$), it is satisfied} $$|Q(H)|=\sum_{k=m+1}^M Q_l(Z_k).\eqno(\star)$$Note that $1-Q_l(Z_0\cap H^c)=Q_l([\psi_{\mathcal Z}(H)]^c)$ implies $|Q(H)|=\sum_{k=1}^M Q_l(Z_k)$ when $H$ is a non-empty event with $H\subseteq Z_0$. Then, owing to relation $(ii)$, condition $(\star)$ is true when $m=0$, while condition $(\star)$ is implied by $(iii)$ through an inductive argument for $m\geq 1$. Moreover, it is immediate to verify conditions $(i)$, $(ii)$, and $(iii)$ for $P_{\psi_{\mathcal Z}}$, and ${\mathcal W}$ is the algebra generated by $ {\mathcal Z}$.
\bigskip
\subsection{ Proof Proposition 14}
Let $Q=P_{\mathcal Y}$ for a suitable not increasing family ${\mathcal Y}$. For every finite subset $I$ of $J$, since $Y_{\Omega\setminus\{\omega_i\}_{i\in I}}\geq Y_{\Omega\setminus\{\omega_i\}}$ and $Y_{\Omega\setminus\{\omega_i\}}=0$ on $\Omega\setminus\{\omega_i\}$ for every $i\in I$, it yields $$Y_{\Omega\setminus\{\omega_i\}_{i\in I}}\geq \sum_{i\in I} Y_{\Omega\setminus\{\omega_i\}}.$$ 
Then $$|Q({\Omega\setminus\{\omega_i\}_{i\in I}})|=E^{Q_l}[Y_{\Omega\setminus\{\omega_i\}_{i\in I}}]\geq \sum_{i\in I} E^{Q_l}[Y_{\Omega\setminus\{\omega_i\}}]= \sum_{i\in I}|Q(\Omega\setminus\{\omega_i\})| $$
Conversely, pose
$$Y_H=r_HI_{H^c}+\sum_{\omega_i\in H^c}Y_{\Omega\setminus\{\omega_i\}}$$
where $Y_{\Omega\setminus\{\omega_i\}}={{|Q(\Omega\setminus\{\omega_i\})|}\over{Q_l(\{\omega_i\})}}I_{\{\omega_i\}}$ and $r_H$ verifies $$r_HQ_l(H^c)=|Q(H)|-\sum_{\omega_i\in H^c}|Q(\Omega\setminus\{\omega_i\})|.$$ Owing to Definition 7, $(Y_H)_H$ is a not increasing family of random variables and it follows $Q=(Q_l)_{\mathcal Y}$. 
\bigskip
\subsection{ Proof Proposition 18}

Let $H$ be an event not  $P_{\psi_{\mathcal{Z}}}$-negligible and $H\subseteq Z_n$
for some $n$. If $A$ is an event contained in some $Z_m$, for $m\neq n$, it holds 
\begin{eqnarray*}
P_{\psi_{\mathcal{Z}}}(A|H) &=& \left[ \frac{P(A\cap[\psi_{\mathcal Z}(H)]^c)}{P([\psi_{\mathcal Z}(A)]^c}, \frac{P([\psi_{\mathcal Z}(A)]^c\cap[\psi_{\mathcal Z}(H)]^c)}{P([\psi_{\mathcal Z}(A)]^c} \right]
\\
&=&\left[ \frac{P(A)}{P(H) + P(H_{\psi_{\mathcal Z}})]^c}, \frac{P([\psi_{\mathcal Z}(A)]^c)-P([\psi_{\mathcal Z}(A)]^c\cap \psi_{\mathcal Z}(H))}{P(H) + P(H_{\psi_{\mathcal Z}})} \right]\\
&=&\left[ \frac{P(A)}{P(H) +1-P(Z_n)}, \frac{P([\psi_{\mathcal Z}(A)]^c)-P( \psi_{\mathcal Z}(H))}{P(H) +1-P(Z_n)} \right]\\
&=&\left[ \frac{P(A)}{P(H) +1-P(Z_n)}, \frac{1-P(A^c\cap Z_m)-(P(Z_n)-P(H))}{P(H) +1-P(Z_n)} \right]\\
\end{eqnarray*}
that is the relation $(22)$ is true.
For $n=m$, it results 
\begin{eqnarray*}
P_{\psi_{\mathcal{Z}}}(A|H) &=& \left[ \frac{P(A\cap H)}{P(H) + 1 - P(Z_n)}, \frac{%
P(A\cap H) + P(A_{\psi_{\mathcal Z}}\cap H_{\psi_{\mathcal Z}})}{P(H) + 1 - P(Z_n)} \right]\\ &=& \left[ \frac{P(A\cap H)}{P(H) + 1 - P(Z_n)}, \frac{%
P(A\cap H) + 1 - P(Z_n)}{P(H) + 1 - P(Z_n)} \right]
\end{eqnarray*}
because $A_{\psi_{\mathcal Z}}\cap H_{\psi_{\mathcal Z}}=Z_n^c$. Proposition $18$ is so proven.

\medskip
\begin{thebibliography}{9}

\bibitem{Birk} Birkhoff, G. and Von Neumann J. (1936): The Logic of Quantum Mechanics, {\it Annals of Mathematics}, {\bf 37}, 823-843. 

\bibitem{Boole} Boole, G. (1958,[1854]): {\it An Investigation of the Laws of Thought on Which are Founded the mathematical Theories of Logic and Probabilities} Reprinted with corrections, Dover Publications, (reissued by Cambridge University Press), (2009).

\bibitem{Dempster} Dempster, A. (1967): Upper and lower probabilities induced by a multivalued mapping, {\it The Annals of Mathematical Statistics}, {\bf 38}, 325-339.

\bibitem{Fagin} Fagin, R., Halpern J.Y., Moses Y. and Vardi, M.Y. (1995): {\it Reasoning about knowledge}. MIT press, Paperback edition (2003).

\bibitem{Itzhah} 
Itzhak, G. and Schmeidler D. (1989): Maxmin expected utility with non-unique prior, {\it Journal of Mathematical Economics}, {\bf 18}, 141-153. 
 
\bibitem{Good}  Good, I. (1982): The measure of a non-measurable set, in { Logic, Methodology and Philosophy of Science} (edited by Nagel, Suppes and Taski). Stanford university Press, 319-329.

\bibitem{Keynes} Keynes, J.M. (1921): {\it A Treatise On Probability}, London, Macmillan and Co.

\bibitem{Koopman} Koopman, B. (1940): The axioms and algebra of intuitive probability. {\it Annals of Mathematics}, {\bf 41}, 269-292.

\bibitem{Kyburg} Kyburg, H. (1963): Probability and randomness, {\it Theoria}, {\bf 29}, 27–55.

\bibitem{Ipcc} IPCC Fifth Assessment Report (AR5): Climate Change 2014: Synthesis Report. Contribution of Working Groups I, II and III to the Fifth Assessment Report of the Intergovernmental Panel on Climate Change.

\bibitem{Ipcc2} IPCC: Climate Change 2023: Synthesis Report. Contribution of Working Groups I, II and III to the Sixth Assessment Report of the Intergovernmental Panel on Climate Change [Core Writing Team, H. Lee and J. Romero (eds.)]. IPCC, Geneva, Switzerland, pp. 35-115, doi: 10.59327/IPCC/AR6-9789291691647

\bibitem{Ramsey} Ramsey, F. (1926):  {\it Truth and Probability}, in Ramsey, 1931, {\it The foundation of mathematics and other Logical Essays}, Ch VII, 156-198. Edited by Brainthwaite B.: Kegan, Paul, Trench, Trubner and Co., Harcourt.

\bibitem{Shaffer}
Shaffer, G. (1976): {\it A mathematical theory of evidence}, vol. 1. Princeton University Press.

\bibitem{Smets} Smets, P. (1981): The degree of belief in a fuzzy event, {\it Information Sciences}, {\bf 25}, 1–19.

\bibitem{Stone} Stone, M. (1938): The representation of Boolean algebras,  {\it Bulletin of the American Mathematical Society}, {\bf 44}, 807–816.

\end{thebibliography}
\end{document}